# Confidence bands for convex median curves using sign-tests

## Lutz Dümbgen[1]

*University of Bern*

**Abstract:** Suppose that one observes pairs $(x_1, Y_1)$, $(x_2, Y_2)$, ..., $(x_n, Y_n)$, where $x_1 \le x_2 \le \cdots \le x_n$ are fixed numbers, and $Y_1, Y_2, \ldots, Y_n$ are independent random variables with unknown distributions. The only assumption is that $\mathrm{Median}(Y_i) = f(x_i)$ for some unknown convex function $f$. We present a confidence band for this regression function $f$ using suitable multiscale sign-tests. While the exact computation of this band requires $O(n^4)$ steps, good approximations can be obtained in $O(n^2)$ steps. In addition the confidence band is shown to have desirable asymptotic properties as the sample size $n$ tends to infinity.

## 1. Introduction

Suppose that we are given data vectors $\mathbf{x}, \mathbf{Y} \in \mathbb{R}^n$, where $\mathbf{x}$ is a fixed vector with components $x_1 \le x_2 \le \cdots \le x_n$, and $\mathbf{Y}$ has independent components $Y_i$ with unknown distributions. We assume that

(1) $$\mathrm{Median}(Y_i) = f(x_i)$$

for some unknown convex function $f : \mathbb{R} \to \bar{\mathbb{R}}$, where $\bar{\mathbb{R}}$ denotes the extended real line $[-\infty, \infty]$. To be precise, we assume that $f(x_i)$ is *some* median of $Y_i$. In what follows we present a confidence band $(\hat{L}, \hat{U})$ for $f$. That means, $\hat{L} = \hat{L}(\cdot \,|\, \mathbf{x}, \mathbf{Y}, \alpha)$ and $\hat{U} = \hat{U}(\cdot \,|\, \mathbf{x}, \mathbf{Y}, \alpha)$ are data-dependent functions from $\mathbb{R}$ into $\bar{\mathbb{R}}$ such that

(2) $$\mathbb{P}\left( \hat{L}(x) \le f(x) \le \hat{U}(x) \text{ for all } x \in \mathbb{R} \right) \ge 1 - \alpha$$

for a given level $\alpha \in (0, 1)$.

Our confidence sets are based on a multiscale sign-test. A similar method has been applied by Dümbgen and Johns [2] to treat the case of *isotonic* regression functions, and the reader is referred to that paper for further references. The remainder of the present paper is organized as follows: Section 2 contains the explicit definition of our sign-test statistic and provides some critical values. A corresponding confidence band $(\hat{L}, \hat{U})$ is described in Section 3. This includes exact algorithms for the computation of the upper bound $\hat{U}$ and the lower bound $\hat{L}$ whose running time is of order $O(n^4)$ and $O(n^3)$, respectively. For large data sets these computational complexities are certainly too high. Therefore we present approximate solutions in Section 4 whose running time is of order $O(n^2)$. In Section 5 we discuss the asymptotic behavior of the width of our confidence band as the sample size $n$ tends to

---

infinity. Finally, in Section 6 we illustrate our methods with simulated and real data.

Explicit computer code (in *MatLab*) for the procedures of the present paper as well as of Dümbgen and Johns [2] may be downloaded from the author's homepage.

## 2. Definition of the test statistic

Given any candidate $g : \mathbb{R} \to \bar{\mathbb{R}}$ for $f$ we consider the sign vectors $\underline{\mathrm{sign}}(\mathbf{Y} - g(\mathbf{x}))$ and $\underline{\mathrm{sign}}(g(\mathbf{x}) - \mathbf{Y})$, where $g(\mathbf{x}) := (g(x_i))_{i=1}^n$ and

$$\underline{\mathrm{sign}}(x) := 1\{x > 0\} - 1\{x \le 0\} \quad \text{for } x \in \bar{\mathbb{R}},$$
$$\underline{\mathrm{sign}}(\mathbf{v}) := \big(\underline{\mathrm{sign}}(v_i)\big)_{i=1}^n \quad \text{for } \mathbf{v} = (v_i)_{i=1}^n \in \bar{\mathbb{R}}^n.$$

This non-symmetric definition of the sign function is necessary in order to deal with possibly non-continuous distributions. Whenever the vector $\underline{\mathrm{sign}}(\mathbf{Y} - g(\mathbf{x}))$ or $\underline{\mathrm{sign}}(g(\mathbf{x}) - \mathbf{Y})$ contains "too many" ones in some region, the function $g$ is rejected. Our confidence set for $f$ comprises all convex functions $g$ which are not rejected.

Precisely, let $T_o : \{-1,1\}^n \to \mathbb{R}$ be some test statistic such that $T_o(\boldsymbol{\sigma}) \le T_o(\tilde{\boldsymbol{\sigma}})$ whenever $\boldsymbol{\sigma} \le \tilde{\boldsymbol{\sigma}}$ component-wise. Then we define

$$T(\mathbf{v}) := \max\big\{T_o(\underline{\mathrm{sign}}(\mathbf{v})), T_o(\underline{\mathrm{sign}}(-\mathbf{v}))\big\}$$

for $v \in \bar{\mathbb{R}}^n$. Let $\boldsymbol{\xi} \in \{-1,1\}^n$ be a Rademacher vector, i.e. a random vector with independent components $\xi_i$ which are uniformly distributed on $\{-1,1\}$. Further let $\kappa = \kappa(n, \alpha)$ be the smallest $(1-\alpha)$–quantile of $T(\boldsymbol{\xi})$. Then

$$\mathbb{P}\big(T(\mathbf{Y} - f(\mathbf{x})) \le \kappa\big) \ge \mathbb{P}(T(\boldsymbol{\xi}) \le \kappa) \ge 1 - \alpha;$$

see Dümbgen and Johns [2]. Consequently the set

$$C(\mathbf{x}, \mathbf{Y}, \alpha) := \big\{\text{convex } g : T(\mathbf{Y} - g(\mathbf{x})) \le \kappa\big\}$$

contains $f$ with probability at least $1 - \alpha$.

As for the test statistic $T_o$, let $\psi$ be the triangular kernel function given by

$$\psi(x) := \max(1 - |x|, 0).$$

Then we define

$$T_o(\boldsymbol{\sigma}) := \max_{d=1,\ldots,\lfloor (n+1)/2 \rfloor}\left(\max_{j=1,\ldots,n} T_{d,j}(\boldsymbol{\sigma}) - \Gamma\Big(\frac{2d-1}{n}\Big)\right),$$

where

$$\Gamma(u) := (2\log(e/u))^{1/2},$$
$$T_{d,j}(\boldsymbol{\sigma}) := \beta_d \sum_{i=1}^n \psi\Big(\frac{i-j}{d}\Big)\sigma_i \quad \text{with } \beta_d := \Big(\sum_{i=1-d}^{d-1} \psi\Big(\frac{i}{d}\Big)^2\Big)^{-1/2}.$$

Note that $T_{d,j}(\boldsymbol{\sigma})$ is measuring whether $(\sigma_i)_{j-d<i<j+d}$ contains suspiciously many ones. Thus $d$ and $j$ can be viewed as scale and location parameter, respectively. The normalizing constant $\beta_d$ is chosen such that the standard deviation of $T_{d,j}(\boldsymbol{\xi})$ is not greater than one, with equality if $d \le j \le n+1-d$. The additive correction term $\Gamma((2d-1)/n)$ is justified by results of Dümbgen and Spokoiny [3] about multiscale testing. In fact, Theorem 6.1 of Dümbgen and Spokoiny [3] and Donsker's invariance principle for partial sums of the Rademacher vector $\boldsymbol{\xi}$ together imply that the distribution of $T(\boldsymbol{\xi})$ converges weakly to a probability distribution on $[0, \infty)$ as $n \to \infty$.

TABLE 1
*Critical values $\kappa(n, \alpha)$*

| $\alpha$ | Sample size $n$ | | | | | | | | |
|---|---|---|---|---|---|---|---|---|---|
| | 100 | 200 | 300 | 500 | 700 | 1000 | 2000 | 5000 | 10000 |
| 0.50 | 0.054 | 0.124 | 0.152 | 0.188 | 0.216 | 0.232 | 0.279 | 0.333 | 0.362 |
| 0.10 | 0.792 | 0.860 | 0.867 | 0.904 | 0.902 | 0.915 | 0.970 | 0.991 | 1.021 |
| 0.05 | 1.035 | 1.102 | 1.102 | 1.135 | 1.136 | 1.152 | 1.229 | 1.231 | 1.246 |

Explicit formulae for quantiles of the limiting distribution of $T(\boldsymbol{\xi})$ are not available. Therefore we list some quantiles of $T(\boldsymbol{\xi})$ for various values of $n$ and $\alpha$ in Table 1. Each quantile has been estimated in 19999 Monte Carlo simulations.

## 3. Definition and exact computation of a band

In principle one could define a confidence band $(\tilde{L}, \tilde{U})$ via

$$\tilde{L} := \inf\{g \in C(\mathbf{x}, \mathbf{Y}, \alpha)\}$$
$$= \inf\left\{\text{convex } g \,:\, T_o(\underline{\text{sign}}(\mathbf{Y} - g(\mathbf{x})) \leq \kappa, T_o(\underline{\text{sign}}(g(\mathbf{x}) - \mathbf{Y}) \leq \kappa\right\},$$
$$\tilde{U} := \sup\{g \in C(\mathbf{x}, \mathbf{Y}, \alpha)\}$$
$$= \sup\left\{\text{convex } g \,:\, T_o(\underline{\text{sign}}(\mathbf{Y} - g(\mathbf{x})) \leq \kappa, T_o(\underline{\text{sign}}(g(\mathbf{x}) - \mathbf{Y}) \leq \kappa\right\}.$$

Throughout this paper maxima or minima of functions are defined pointwise. Unfortunately, the explicit computation of $(\tilde{L}, \tilde{U})$ is far from trivial. Therefore we modify the latter definition and compute a band $(\hat{L}, \hat{U})$ in two steps. Our upper boundary is given by

$$\hat{U} := \max\left\{\text{convex } g \,:\, T_o(\underline{\text{sign}}(g(\mathbf{x}) - \mathbf{Y}) \leq \kappa\right\}.$$

Thus we just drop the constraint $T_o(\underline{\text{sign}}(\mathbf{Y} - g(\mathbf{x}))) \leq \kappa$ in the definition of $\tilde{U}$ and obtain $\hat{U} \geq \tilde{U}$. With $\hat{U}$ at hand, our lower boundary is defined as

$$\hat{L} := \min\left\{\text{convex } g \,:\, g \leq \hat{U}, T_o(\underline{\text{sign}}(\mathbf{Y} - g(\mathbf{x}))) \leq \kappa\right\}.$$

Here we replace the constraint $T_o(\underline{\text{sign}}(g(\mathbf{x}) - \mathbf{Y})) \leq \kappa$ in the definition of $\tilde{L}$ with the weaker constraint $g \leq \hat{U}$ and obtain $\hat{L} \leq \tilde{L}$. In what follows we concentrate on the computation of the corresponding vectors $\hat{\boldsymbol{L}} = (\hat{L}_i)_{i=1}^\infty = \hat{L}(\mathbf{x})$ and $\hat{\boldsymbol{U}} = (\hat{U}_i)_{i=1}^n = \hat{U}(\mathbf{x})$.

### 3.1. Computation of $\hat{U}$

**A simplified expression for $\hat{U}$.** To determine $\hat{U}$ it suffices to consider the class $\mathcal{G}$ consisting of the following convex functions $g_{j,k}$: For $1 \leq j < k \leq n$ with $x_j < x_k$ define

$$g_{j,k}(x) := Y_j + \frac{Y_k - Y_j}{x_k - x_j}(x - x_j),$$

describing a straight line connecting the data points $(x_j, Y_j)$ and $(x_k, Y_k)$. Moreover, for $j, k \in \{1, \ldots, n\}$ let

$$g_{0,k}(x) := \begin{cases} \infty & \text{if } x < x_k, \\ Y_k & \text{if } x = x_k, \\ -\infty & \text{if } x > x_k. \end{cases}$$

$$g_{j,n+1}(x) := \begin{cases} -\infty & \text{if } x < x_j, \\ Y_j & \text{if } x = x_j, \\ \infty & \text{if } x > x_j. \end{cases}$$

Then

(3) $$\hat{U} = \max\{g : g \in \mathcal{G}, T_o(\underline{\text{sign}}(g(\mathbf{x}) - \mathbf{Y})) \leq \kappa\}.$$

For let $g$ be any convex function such that $T_o(\underline{\text{sign}}(g(\mathbf{x}) - \mathbf{Y})) \leq \kappa$. Let $\tilde{g}$ be the largest convex function such that $\tilde{g}(x_i) \leq Y_i$ for all indices $i$ with $g(x_i) \leq Y_i$. This function $\tilde{g}$ is closely related to the convex hull of all data points $(x_i, Y_i)$ with $g(x_i) \leq Y_i$. Obviously, $\tilde{g} \geq g$ and $T_o(\underline{\text{sign}}(\tilde{g}(\mathbf{x}) - \mathbf{Y})) = T_o(\underline{\text{sign}}(g(\mathbf{x}) - \mathbf{Y}))$. Let $\omega(1) < \cdots < \omega(m)$ be indices such that $x_{\omega(1)} < \cdots < x_{\omega(m)}$ and

$$\{(x, \tilde{g}(x)) : x \in \mathbb{R}\} \cap \{(x_i, Y_i) : 1 \leq i \leq n\} = \{(x_{\omega(\ell)}, Y_{\omega(\ell)}) : 1 \leq \ell \leq m\}.$$

With $\omega(0) := 0$ and $\omega(m+1) := n+1$ one may write $\tilde{g}$ as the maximum of the functions $g_{\omega(\ell-1),\omega(\ell)}$, $1 \leq \ell \leq m+1$, all of which satisfy the inequality $T_o(\underline{\text{sign}}(g_{\omega(\ell-1),\omega(\ell)}(\mathbf{x}) - \mathbf{Y})) \leq T_o(\underline{\text{sign}}(\tilde{g}(\mathbf{x}) - \mathbf{Y})) \leq \kappa$. Figure 1 illustrates these considerations.

**Computational complexity.** As we shall explain in Section 3.3, the computation of $T_o(\underline{\text{sign}}(g(\mathbf{x}) - \mathbf{Y}))$ for one single candidate function $g \in \mathcal{G}$ requires $O(n^2)$ steps. In case of $T_o(\underline{\text{sign}}(g(\mathbf{x}) - \mathbf{Y})) \leq \kappa$ we have to replace $\hat{U}$ with the vector $\left(\max(g(x_i), \hat{U}_i)\right)_{i=1}^n$ in another $O(n)$ steps. Consequently, since $\mathcal{G}$ contains at most $n(n-1)/2 + 2n = O(n^2)$ functions, the computation of $\hat{U}$ requires $O(n^4)$ steps.

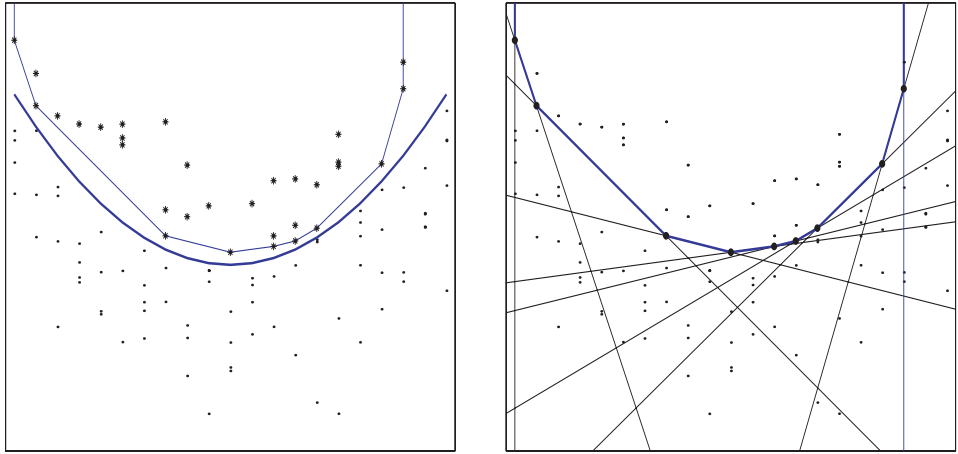

FIG 1. *A function $g$ and its associated function $\tilde{g}$.*

### 3.2. Computation of $\hat{L}$

From now on we assume that $\hat{U}$ is nontrivial, i.e. that $\hat{U}_i = \hat{U}(x_i) < \infty$ for some value $x_i$. Moreover, letting $x_{\min}$ and $x_{\max}$ be the smallest and largest such value, we assume that $x_{\min} < x_{\max}$. Finally let $T_o(\underline{\text{sign}}(\mathbf{Y} - \hat{U}(\mathbf{x}))) \leq \kappa$. Otherwise the confidence set $C(\mathbf{x}, \mathbf{Y}, \alpha)$ would be empty, meaning that convexity of the median function is not plausible.

**Simplified formulae for $\hat{L}$.**

Similarly as in the previous section, one may replace the set of all convex functions with a finite subset $\mathcal{H} = \mathcal{H}(\hat{U})$. First of all let $h$ be any convex function such that $h \leq \hat{U}$ and $T_o(\underline{\text{sign}}(\mathbf{Y} - h(\mathbf{x}))) \leq \kappa$. For any real number $t$ let $z := h(t)$. Now let $\tilde{h} = \tilde{h}_{t,z}$ be the largest convex function such that $\tilde{h} \leq \hat{U}$ and $\tilde{h}(t) = z$. Obviously $\tilde{h} \geq h$, whence $T_o(\underline{\text{sign}}(\mathbf{Y} - \tilde{h}(\mathbf{x}))) \leq \kappa$. Consequently,

$$\hat{L}(t) = \inf\left\{z \in \mathbb{R} : T_o(\underline{\text{sign}}(\mathbf{Y} - \tilde{h}_{t,z}(\mathbf{x}))) \leq \kappa\right\}. \tag{4}$$

Figure 2 illustrates the definition of $\tilde{h}_{t,z}$. Note that $\tilde{h}_{t,z}$ is given by the convex hull of the point $(t, z)$ and the epigraph of $\hat{U}$, i.e. the set of all pairs $(x, y) \in \mathbb{R}^2$ such that $\hat{U}(x) \leq y$.

Starting from equation (4) we derive a computable expression for $\hat{L}$. For that purpose we define tangent parameters as follows: Let $\mathcal{J}$ be the set of all indices $j \in \{1, \ldots, n\}$ such that $\hat{U}(x_j) \geq Y_j$. For $j \in \mathcal{J}$ define

$$s_j^l := \begin{cases} -\infty & \text{if } x_j \leq x_{\min}, \\ \max_{x_i < x_j} \dfrac{Y_j - \hat{U}(x_i)}{x_j - x_i} & \text{else,} \end{cases}$$

$$a_j^l := \begin{cases} x_j & \text{if } x_j \leq x_{\min}, \\ \arg\max_{x_i < x_j} \dfrac{Y_j - \hat{U}(x_i)}{x_j - x_i} & \text{else,} \end{cases}$$

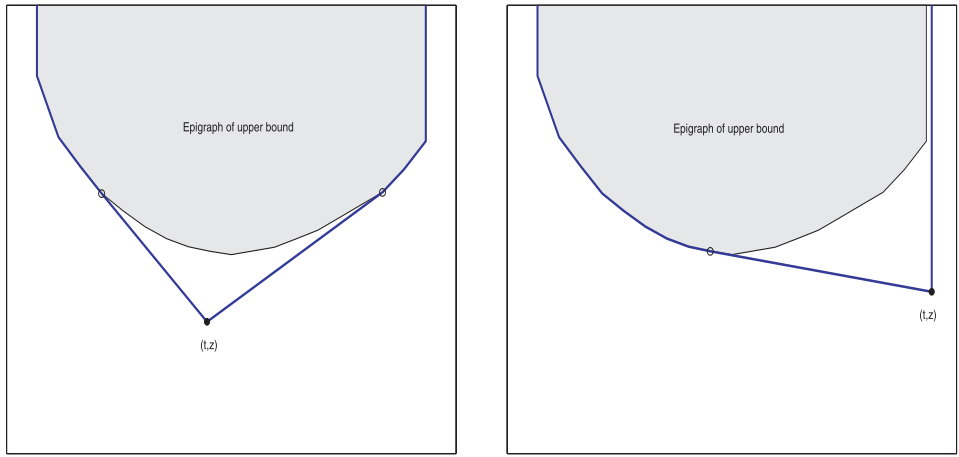

FIG 2. *The extremal function $\tilde{h}_{t,z}$ of two points $(t, z)$.*

$$s_j^{\mathrm{r}} := \begin{cases} \infty & \text{if } x_j \geq x_{\max}, \\ \min_{x_k > x_j} \dfrac{\hat{U}(x_k) - Y_j}{x_k - x_j} & \text{else,} \end{cases}$$

$$a_j^{\mathrm{r}} := \begin{cases} x_j & \text{if } x_j \geq x_{\max}, \\ \arg\min_{x_k > x_j} \dfrac{\hat{U}(x_k) - Y_j}{x_k - x_j} & \text{else.} \end{cases}$$

With these parameters we define auxiliary tangent functions

$$h_j^{\mathrm{l}}(x) := \begin{cases} \hat{U}(x) & \text{if } x < a_j^{\mathrm{l}}, \\ Y_j + s_j^{\mathrm{l}}(x - x_j) & \text{if } x \geq a_j^{\mathrm{l}}, \end{cases}$$

$$h_j^{\mathrm{r}}(x) := \begin{cases} Y_j + s_j^{\mathrm{r}}(x - x_j) & \text{if } x \leq a_j^{\mathrm{r}}, \\ \hat{U}(x) & \text{if } x > a_j^{\mathrm{r}}. \end{cases}$$

Figure 3 depicts these functions $h_j^{\mathrm{l}}$ and $h_j^{\mathrm{r}}$. Note that

$$h_j^{\mathrm{l}}(x) = \begin{cases} \max\left\{h(x) : h \text{ convex}, h \leq \hat{U}, h(x_j) \leq Y_j\right\} & \text{if } x \leq x_j, \\ \min\left\{h(x) : h \text{ convex}, h \leq \hat{U}, h(x_j) \geq Y_j\right\} & \text{if } x \geq x_j, \end{cases}$$

$$h_j^{\mathrm{r}}(x) = \begin{cases} \min\left\{h(x) : h \text{ convex}, h \leq \hat{U}, h(x_j) \geq Y_j\right\} & \text{if } x \leq x_j, \\ \max\left\{h(x) : h \text{ convex}, h \leq \hat{U}, h(x_j) \leq Y_j\right\} & \text{if } x \geq x_j, \end{cases}$$

In particular, $h_j^{\mathrm{l}}(x_j) = h_j^{\mathrm{r}}(x_j) = Y_j$. In addition we define $h_0^{\mathrm{l}}(x) := h_{n+1}^{\mathrm{r}}(x) := -\infty$. Then we set

$$h_{j,k} := \max(h_j^{\mathrm{l}}, h_k^{\mathrm{r}}) \quad \text{and} \quad \mathcal{H} := \{h_{j,k} : j \in \{0\} \cup \mathcal{J}, k \in \mathcal{J} \cup \{n+1\}\}.$$

This class $\mathcal{H}$ consists of at most $(n+1)^2$ functions, and elementary considerations show that

(5) $$\hat{L} = \min\left\{h \in \mathcal{H} : T_o(\underline{\mathrm{sign}}(\mathbf{Y} - h(\mathbf{x}))) \leq \kappa\right\}.$$

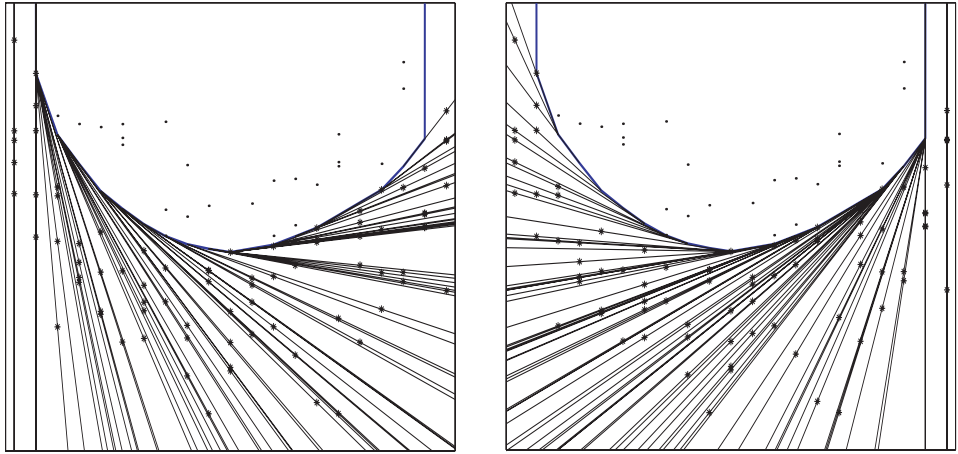

FIG 3. *The tangent functions $h_j^{\mathrm{l}}$ and $h_k^{\mathrm{r}}$.*

**Computational complexity.** Note first that any pair $(a_j^{\mathrm{w}}, s_j^{\mathrm{w}})$ may be computed in $O(n)$ steps. Consequently, before starting with $\hat{\boldsymbol{L}}$ we may compute all tangent parameters in time $O(n^2)$. Then Equation (5) implies that $\hat{\boldsymbol{L}}$ may be computed in $O(n^4)$ steps. However, this can be improved considerably. The reason is, roughly saying, that for fixed $j$, one can determine the smallest function $h_k^{\mathrm{r}}$ such that $T_o(\mathrm{sign}(\mathbf{Y} - h_{j,k}(\mathbf{x}))) \le \kappa$ in $O(n^2)$ steps, as explained in the subsequent section. Hence a proper implementation lets us compute $\hat{\boldsymbol{L}}$ in $O(n^3)$ steps.

### 3.3. An auxiliary routine

In this section we show that the value of $T_o(\boldsymbol{\sigma})$ can be computed in $O(n^2)$ steps. More generally, we consider $n$–dimensional sign vectors $\boldsymbol{\sigma}^{(0)}, \boldsymbol{\sigma}^{(1)}, \ldots, \boldsymbol{\sigma}^{(q)}$ such that for $1 \le \ell \le q$ the vectors $\boldsymbol{\sigma}^{(\ell-1)}$ and $\boldsymbol{\sigma}^{(\ell)}$ differ exactly in one component, say,

$$\sigma_{\omega(\ell)}^{(\ell-1)} = 1 \quad \text{and} \quad \sigma_{\omega(\ell)}^{(\ell)} = -1$$

for some index $\omega(\ell) \in \{1, \ldots, n\}$. Thus $\boldsymbol{\sigma}^{(0)} \ge \boldsymbol{\sigma}^{(1)} \ge \cdots \ge \boldsymbol{\sigma}^{(q)}$ component-wise. In particular, $T_o(\boldsymbol{\sigma}^{(\ell)})$ is non-increasing in $\ell$. It is possible to determine the number

$$\ell_* := \min\left(\left\{\ell \in \{0, \ldots, q\} : T_o(\boldsymbol{\sigma}^{(\ell)}) \le \kappa\right\} \cup \{\infty\}\right)$$

in $O(n^2)$ steps as follows:

**Algorithm.** We use three vector variables $\boldsymbol{S}$, $\boldsymbol{S}^{(0)}$ and $\boldsymbol{S}^{(1)}$ plus two scalar variables $\ell$ and $d$. While running the algorithm the variable $\boldsymbol{S}$ contains the current vector $\boldsymbol{\sigma}^{(\ell)}$, while

$$\boldsymbol{S}^{(0)} = \Big(\sum_{i \in [j-d+1, j+d-1]} \boldsymbol{S}_i\Big)_{j=1}^n,$$

$$\boldsymbol{S}^{(1)} = \Big(\sum_{i \in [j-d+1, j+d-1]} (d - |j-i|) S_i\Big)_{j=1}^n.$$

**Initialisation.**

$$\ell \leftarrow 0, \quad d \leftarrow 1 \quad \text{and} \quad \left.\begin{matrix} \boldsymbol{S} \\ \boldsymbol{S}^{(0)} \\ \boldsymbol{S}^{(1)} \end{matrix}\right\} \leftarrow \boldsymbol{\sigma}^{(0)}.$$

**Induction step.** Check whether

(6) $$\max_{i=1,\ldots,n} S_i^{(1)} \le \Big(\sum_{i=1-d}^{d-1}(d-i)^2\Big)^{1/2}\big(\Gamma((2d-1)/n) + \kappa\big)$$
$$= ((2d^2 + 1)d/3)^{1/2}\big(\Gamma((2d-1)/n) + \kappa\big).$$

- If (6) is fulfilled and $d < \lfloor(n+1)/2\rfloor$, then

$$d \leftarrow d + 1,$$

$$S_i^{(0)} \leftarrow \begin{cases} S_i^{(0)} + S_{i+d-1} & \text{for } i < d, \\ S_i^{(0)} + S_{i+1-d} + S_{i+d-1} & \text{for } d \le i \le n+1-d, \\ S_i^{(0)} + S_{i+1-d} & \text{for } i > n+1-d, \end{cases}$$

$$\boldsymbol{S}^{(1)} \leftarrow \boldsymbol{S}^{(1)} + \boldsymbol{S}^{(0)}.$$

- If (6) is fulfilled and $d = \lfloor (n+1)/2 \rfloor$, then

$$\ell_* \leftarrow \ell.$$

- If (6) is violated and $\ell < q$, then

$$\begin{aligned}
\ell &\leftarrow \ell + 1, \\
S_{\omega(\ell)} &\leftarrow -1, \\
S_i^{(0)} &\leftarrow S_i^{(0)} - 2 \quad \text{and} \\
S_i^{(1)} &\leftarrow S_i^{(1)} - 2(d - |i - \omega(\ell)|) \quad \text{for } \omega(\ell) - d < i < \omega(\ell) + d.
\end{aligned}$$

- If Condition (6) is violated but $\ell = q$, then $T_o(\boldsymbol{\sigma}^{(q)}) > \kappa$, and

$$\ell_* \leftarrow \infty.$$

As for the running time of this algorithm, note that each induction step requires $O(n)$ operations. Since either $d$ or $\ell$ increases each time by one, the algorithm terminates after at most $n + q + 1 \leq 2n + 1$ induction steps. Together with $O(n)$ operations for the initialisation we end up with total running time $O(n^2)$.

## 4. Approximate solutions

**Approximation of $\hat{U}$.** Recall that the exact computation of $\hat{U}$ involves testing whether a straight line given by a function $g(\cdot)$ and touching one or two data points $(x_i, Y_i)$ satisfies the inequality $T_o(\underline{\text{sign}}(g(\mathbf{x}) - \mathbf{Y})) \leq \kappa$. The idea of our approximation is to restrict our attention to straight lines whose slope belongs to a given finite set.

**Step 1.** At first we consider the straight lines $g_{0,k}$ instroduced in section 3.1, all having slope $-\infty$. Let $\omega(1), \ldots, \omega(n)$ be a list of $\{1, \ldots, n\}$ such that $g_{0,\omega(1)} \leq \cdots \leq g_{0,\omega(n)}$. In other words, for $1 < \ell \leq n$ either $x_{\omega(\ell-1)} < x_{\omega(\ell)}$, or $x_{\omega(\ell-1)} = x_{\omega(\ell)}$ and $Y_{\omega(\ell-1)} \leq Y_{\omega(\ell)}$. With the auxiliary procedure of Section 3.3 we can determine the the smallest number $\ell_*$ such that $T_o(\underline{\text{sign}}(g_{0,\omega(\ell_*)}(\mathbf{x}) - \mathbf{Y})) \leq \kappa$ in $O(n^2)$ steps. We write $G_0 := g_{0,\omega(\ell_*)}$. Note that $x_{\omega(\ell_*)}$ is equal to $x_{\min} = \min\{x : \hat{U}(x) < \infty\}$.

**Step 2.** For any given slope $s \in \mathbb{R}$ let $a(s)$ be the largest real number such that the sign vector

$$\boldsymbol{\sigma}(s) := \big(\underline{\text{sign}}(Y_i - a(s) - sx_i)\big)_{i=1}^n$$

satisfies the inequality $T_o(\boldsymbol{\sigma}(s)) \leq \kappa$. This number can also be determined in time $O(n^2)$. This time we have to generate and use a list $\omega(1), \ldots, \omega(n)$ of $\{1, 2, \ldots, n\}$ such that $Y_{\omega(\ell)} - sx_{\omega(\ell)}$ is non-increasing in $\ell$.

Now we determine the numbers $a(s_1), \ldots, a(s_{M-1})$ for given slopes $s_1 < \cdots < s_{M-1}$. Then we define

$$G_\ell(x) := a(s_\ell) + s_\ell x \quad \text{for } 1 \leq \ell < M.$$

**Step 3.** Finally we determine the largest function $G_M$ among the degenerate linear functions $g_{1,n+1}, \ldots, g_{n,n+1}$ such that $T_o(\underline{\text{sign}}(G_M(\mathbf{x}) - \mathbf{Y})) \leq \kappa$. This is analogous to Step 1 and yields the number $x_{\max} = \max\{x : \hat{U}(x) < \infty\}$.

**Step 4.** By means of this list of finitely many straight lines $G_0, G_1, \ldots, G_M$ one obtains the lower bound $\hat{U}_* := \max(G_0, G_1, \ldots, G_M)$ for $\hat{U}$. In fact, one could even replace $G_\ell$ with the largest convex function $\tilde{G}_\ell$ such that $\tilde{G}_\ell(x_i) \leq Y_i$ whenever $G_\ell(x_i) \leq Y_i$. Each of these functions can be computed via a suitable variant of the pool-adjacent-violators algorithm in $O(n)$ steps; see Robertson et al. [6].

**Step 5.** To obtain an upper bound $\hat{U}^*$ for $\hat{U}$, for $1 \leq \ell \leq M$ let $H_\ell$ be the smallest *concave* function such that $H_\ell(x_i) \geq Y_i$ whenever $\max(G_{\ell-1}(x_i), G_\ell(x_i)) \geq Y_i$. Again $H_\ell$ may be determined via the pool-adjacent-violators algorithm. Then elementary considerations show that

$$\hat{U} \leq \hat{U}^* := \max\left(\hat{U}_*, H_1, H_2, \ldots, H_M\right).$$

All in all, these five steps require $O(Mn^2)$ steps. By visual inpection of these two curves $\hat{U}_*$ and $\hat{U}^*$ one may opt for a refined grid of slopes or use $\hat{U}^*$ as a surrogate for $\hat{U}$.

**Approximation of $\hat{L}$.** Recall that the exact computation amounts to fixing any function $h_j^l$ und finding the smallest function $h_k^r$ such that $T_o(\underline{\text{sign}}(\mathbf{Y} - h_{j,k}(\mathbf{x}))) \leq \kappa$. Now approximations may be obtained by picking only a subset of the potential indices $j$. In addition, one may fix some functions $h_k^r$ and look for the smallest $h_j^l$ satisfying the constraint $T_o(\underline{\text{sign}}(\mathbf{Y} - h_{j,k}(\mathbf{x}))) \leq \kappa$. Again this leads to approximations $\hat{L}_*$ and $\hat{L}^*$ for $\hat{L}$ such that $\hat{L}_* \leq \hat{L} \leq \hat{L}^*$.

## 5. Asymptotic properties

In this section we consider a triangular array of observations $x_i = x_{n,i}$ and $Y_i = Y_{n,i}$. Our confidence band $(\hat{L}, \hat{U})$ will be shown to have certain consistency properties, provided that $f$ satisfies some smoothness condition, and that the following two requirements are met for some constants $-\infty < a < b < \infty$:

(A1) Let $M_n$ denote the empirical distribution of the design points $x_{n,i}$. That means, $M_n(B) := n^{-1}\#\{i : x_{n,i} \in B\}$ for $B \subset \mathbb{R}$. There is a constant $c > 0$ such that

$$\liminf_{n \to \infty} \frac{M_n[a_n, b_n]}{b_n - a_n} \geq c$$

whenever $a \leq a_n < b_n \leq b$ and $\liminf_{n \to \infty} \log(b_n - a_n)/\log n > -1$.

(A2) All variables $Y_i = Y_{n,i}$ with $x_{n,i} \in [a,b]$ satisfy the following inequalities:

$$\left.\begin{array}{l} \mathbb{P}(Y_i < \mu_i + r) \\ \mathbb{P}(Y_i > \mu_i - r) \end{array}\right\} \geq \frac{1 + H(r)}{2} \quad \text{for any } r > 0,$$

where $H$ is some fixed function on $[0, \infty]$ such that

$$\lim_{r \to 0+} \frac{H(r)}{r} > 0.$$

These conditions (A1) and (A2) are satisfied in various standard models, as pointed out by Dümbgen and Johns [2].

**Theorem 1.** *Suppose that assumptions (A1) and (A2) hold.*

*(a) Let $f$ be linear on $[a, b]$. Then for arbitrary $a < a' < b' < b$,*

$$\left.\begin{array}{l} \sup_{x \in [a,b]} \left(f(x) - \hat{L}(x)\right)^+ \\ \sup_{x \in [a',b']} \left(\hat{U}(x) - f(x)\right)^+ \end{array}\right\} = O_p(n^{-1/2}).$$

(b) Let $f$ be Hölder continuous on $[a,b]$ with exponent $\beta \in (1, 2]$. That means, $f$ is differentiable on $[a,b]$ such that for some constant $L > 0$ and arbitrary $x, y \in [a, b]$,

$$|f'(x) - f'(y)| \le L|x - y|^{\beta-1}.$$

Then for $\rho_n := \log(n+1)/n$ and $\delta_n := \rho_n^{1/(2\beta+1)}$,

$$\left.\begin{array}{c} \sup_{x \in [a,b]} \bigl(f(x) - \hat L(x)\bigr)^+ \\ \sup_{x \in [a+\delta_n, b-\delta_n]} \bigl(\hat U(x) - f(x)\bigr)^+ \end{array}\right\} = O_p\bigl(\rho_n^{\beta/(2\beta+1)}\bigr).$$

Part (a) of this theorem explains the empirical findings in Section 6 that the band $(\hat L, \hat U)$ performs particularly well in regions where the regression function $f$ is linear.

*Proof of Theorem 1, step I.* At first we prove the assertions about $\hat U$. Note that for arbitrary $t, z \in \mathbb{R}$ with $z \le \hat U(t)$ there exist parameters $\mu, \nu \in \mathbb{R}$ such that $z = \mu + \nu t$ and

$$S_{d,j}(\mu, \nu) := \sum_{i=1}^n \psi\Bigl(\frac{i-j}{d}\Bigr) \underline{\mathrm{sign}}(\mu + \nu x_i - Y_i)$$

(7) $$\le \beta_d^{-1}\Bigl(\Gamma\Bigl(\frac{2d-1}{n}\Bigr) + \kappa\Bigr) \quad \text{for any } (d, j) \in \mathcal{T}_n;$$

here $\mathcal{T}_n$ denotes the set of all pairs $(d, j)$ of integers $d > 0$, $j \in [d, n+1-d]$. Therefore it is crucial to have good simultaneous upper bounds for $\bigl|S_{d,j}(\mu, \nu) - \Sigma_{d,j}(\mu, \nu)\bigr|$, where

$$\Sigma_{d,j}(\mu, \nu) := \mathbb{E} S_{d,j}(\mu, \nu) = \sum_{i=1}^n \psi\Bigl(\frac{2d-1}{n}\Bigr)\bigl(2\mathbb{P}(Y_i < \mu + \nu x_i) - 1\bigr).$$

One may write $S_{d,j}(\mu, \nu) = \int g_{d,j,\mu,\nu} \, d\Psi_n$ with the random measure

$$\Psi_n := \sum_{i=1}^n \delta_i \otimes \delta_{x_i} \otimes \delta_{Y_i}$$

and the function

$$(i, x, y) \mapsto g_{d,j,\mu,\nu}(i, x, y) := \psi\Bigl(\frac{i-j}{d}\Bigr) \underline{\mathrm{sign}}(\mu + \nu x - y) \in [-1, 1]$$

on $\mathbb{R}^3$. The family of all these functions $g_{d,j,\mu,\nu}$ is easily shown to be a Vapnik-Cervonenkis subgraph class in the sense of van der Vaart and Wellner [7]. Moreover, $\Psi_n$ is a sum of $n$ stochastically independent random probability measures. Thus well-known results from empirical process theory (cf. Pollard [5]) imply that for arbitrary $\eta \ge 0$,

$$\mathbb{P}\left(\sup_{(d,j) \in \mathcal{T}_n, \, \mu, \nu \in \mathbb{R}} \bigl|S_{d,j}(\mu, \nu) - \Sigma_{d,j}(\mu, \nu)\bigr| \ge n^{1/2}\eta\right)$$

(8) $$\le C \exp(-\eta^2/C),$$

$$\mathbb{P}\left(\sup_{\mu, \nu \in \mathbb{R}} \bigl|S_{d,j}(\mu, \nu) - \Sigma_{d,j}(Y, \mu, \nu)\bigr| \ge d^{1/2}\eta \text{ for some } (d, j) \in \mathcal{T}_n\right)$$

(9) $$\le C \exp(2 \log n - \eta^2/C),$$

where $C \geq 1$ is a universal constant. Consequently, for any fixed $\alpha' > 0$ there is a constant $\tilde{C} > 0$ such that the following inequalities are satisfied simultaneously for arbitrary $(d, j) \in \mathcal{T}_n$ and $(\mu, \nu) \in \mathbb{R}^2$ with probability at least $1 - \alpha'$:

$$
(10) \qquad |S_{d,j}(\mu, \nu) - \Sigma_{d,j}(\mu, \nu)| \leq \begin{cases} \tilde{C} n^{1/2}, \\ \tilde{C} d^{1/2} \log(n+1)^{1/2}. \end{cases}
$$

In what follows we assume (10) for some fixed $\tilde{C}$.

*Proof of part* (a) *for* $\hat{U}$. Suppose that $f$ is linear on $[a, b]$, and let $[a', b'] \subset (a, b)$. By convexity of $\hat{U}$, the maximum of $\hat{U} - f$ over $[a', b']$ is attained at $a'$ or $b'$. We consider the first case and assume that $\hat{U}(a') \geq f(a') + \epsilon_n$ for some $\epsilon_n > 0$. Then there exist $\mu, \nu \in \mathbb{R}$ satisfying (7) such that $\mu + \nu a' = f(a') + \epsilon_n$ and $\nu \leq f'(a')$. In particular, $\mu + \nu x - f(x) \geq \epsilon_n$ for all $x \in [a, a']$. Now we pick a pair $(d_n, j_n) \in \mathcal{T}_n$ with $d_n$ as large as possible such that

$$
[x_{j_n - d_n + 1}, x_{j_n + d_n - 1}] \subset [a, a'].
$$

Assumption (A1) implies that $d_n \geq (c/2 + o(1))n$. Now

$$
\Sigma_{d_n, j_n}(\mu, \nu) \geq \sum_{i = j_n - d_n + 1}^{j_n + d_n - 1} \psi\left(\frac{i - j_n}{d_n}\right) H(\epsilon_n) = d_n H(\epsilon_n)
$$

by assumption (A2). Combining this inequality with (7) and (10) yields

$$
(11) \qquad \beta_{d_n}^{-1}\left(\Gamma\left(\frac{2d_n - 1}{n}\right) + \kappa\right) \geq d_n H(\epsilon_n) - \tilde{C} n^{1/2}.
$$

But $\beta_d^{-1} = 3^{-1/2}(2d - 1)^{1/2} + O(d^{-1/2})$, and $x \mapsto x^{1/2}\Gamma(x)$ is non-decreasing on $(0, 1]$. Hence (11) implies that

$$
\begin{aligned}
H(\epsilon_n) &\leq d_n^{-1}\left((3^{-1/2} + o(1))(2d_n - 1)^{1/2}\left(\Gamma\left(\frac{2d_n - 1}{n}\right) + \kappa\right) + \tilde{C} n^{1/2}\right) \\
&\leq d_n^{-1} n^{1/2}\left((3^{-1/2} + o(1))(\Gamma(1) + \kappa) + \tilde{C}\right) \\
&= O(n^{-1/2}).
\end{aligned}
$$

Consequently, $\epsilon_n = O(n^{-1/2})$.

*Proof of part* (b) *for* $\hat{U}$. Now suppose that $f'$ is Hölder-continuous on $[a, b]$ with exponent $\beta - 1 \in (0, 1]$ and constant $L > 0$. Let $\hat{U}(x) \geq f(x) + \epsilon_n$ for some $x \in [a + \delta_n, b - \delta_n]$ and $\epsilon_n > 0$. Then there are numbers $\mu, \nu \in \mathbb{R}$ satisfying (7) such that $\mu + \nu x = f(x) + \epsilon_n$. Let $(d_n, j_n) \in \mathcal{T}_n$ with $d_n$ as large as possible such that either

$$
f'(x) \leq \nu \quad \text{and} \quad [x_{j_n - d_n + 1}, x_{j_n + d_n - 1}] \subset [x, x + \delta_n],
$$

or

$$
f'(x) \geq \nu \quad \text{and} \quad [x_{j_n - d_n + 1}, x_{j_n + d_n - 1}] \subset [x - \delta_n, x].
$$

Assumption (A1) implies that

$$
d_n \geq (c/2 + o(1))\delta_n n.
$$

Moreover, for any $i \in \{j_n - d_n + 1, j_n + d_n - 1\}$,

$$\begin{aligned}
\mu + \nu x_i - f(x_i) &= \epsilon_n + \int_x^{x_i} (\nu - f'(t))\,dt \\
&\geq \epsilon_n + \int_x^{x_i} (f'(x) - f'(t))\,dt \\
&\geq \epsilon_n - L\int_0^{\delta_n} s^{\beta-1}\,ds \\
&= \epsilon_n - O(\delta_n^\beta),
\end{aligned}$$

so that

$$\Sigma_{d_n,j_n}(\mu,\nu) \geq d_n H\big((\epsilon_n - O(\delta_n^\beta))^+\big).$$

Combining this inequality with (7) and (10) yields

$$\begin{aligned}
H\big((\epsilon_n - O(\delta_n^\beta))^+\big) &\leq d_n^{-1}\beta_{d_n}^{-1}\Big(\Gamma\Big(\frac{2d_n-1}{n}\Big)+\kappa\Big) + \tilde{C}d_n^{-1/2}\log(n+1)^{1/2} \\
&\leq d_n^{-1}\beta_{d_n}^{-1}(2^{1/2}\log(n)^{1/2}+\kappa) + \tilde{C}d_n^{-1/2}\log(n+1)^{1/2} \\
&= O(\delta_n^\beta).
\end{aligned}$$
(12)

This entails that $\epsilon_n$ has to be of order $O(\delta_n^\beta) = O\big(\rho_n^{\beta/(2\beta+1)}\big)$. □

*Proof of Theorem 1, step II.* Now we turn our attention to $\hat{L}$. For that purpose we change the definition of $S_{d,j}(\cdot,\cdot)$ and $\Sigma_{d,j}(\cdot,\cdot)$ as follows: Let $U_n$ be a fixed convex function to be specified later. Then for $(t,z) \in \mathbb{R}^2$ we define $h_{n,t,z}$ to be the largest convex function $h$ such that $h \leq U_n$ and $h(t) \leq z$. This definition is similar to the definition of $\tilde{h}_{t,z}$ in Section 3.2. Indeed, if $\hat{U} \leq U_n$ and $\hat{L}(t) \leq z$, then

$$S_{d,j}(t,z) := \sum_{i=1}^n \psi\Big(\frac{i-j}{d}\Big)\underline{\operatorname{sign}}(Y_i - h_{n,t,z}(x_i))$$

$$\leq \beta_d^{-1}\Big(\Gamma\Big(\frac{2d-1}{n}\Big)+\kappa\Big) \quad \text{for any } (d,j) \in \mathcal{T}_n.$$
(13)

Here we set

$$\Sigma_{d,j}(t,z) := \mathbb{E}S_{d,j}(t,z) = \sum_{i=1}^n \psi\Big(\frac{2d-1}{n}\Big)\big(2\mathbb{P}(Y_i > h_{n,t,z}(x_i))-1\big).$$

Again we may and do assume that (10) is true for some constant $\tilde{C}$.

*Proof of part* (a) *for* $\hat{L}$. Suppose that $f$ is linear on $[a,b]$. We define $U_n(x) := f(x) + \gamma n^{-1/2} + 1\{x \notin [a',b']\}\infty$ with constants $\gamma > 0$ and $a < a' < b' < b$. Since $\liminf_{n\to\infty} \mathbb{P}(\hat{U} \leq U_n)$ tends to one as $\gamma \to \infty$, we may assume that $\hat{U} \leq U_n$. Suppose that $\hat{L}(t) \leq z := f(t) - 2\epsilon_n$ for some $t \in [a,b]$ and $\epsilon_n \geq \gamma n^{-1/2}$. A simple geometrical consideration shows that $h_{n,t,z} \leq f - \epsilon_n$ on an interval $[a'',b''] \subset [a,b]$ of length $b''-a'' \geq (b'-a')/3$. If we pick $(d_n,j_n) \in \mathcal{T}_n$ with $d_n$ as large as possible such that $[x_{j_n-d_n+1}, x_{j_n+d_n-1}] \subset [a'',b'']$, then $d_n \geq (c(b'-a')/6 + o(1))n$. Moreover, (13) and (10) entail (11), whence $\epsilon_n = O(n^{-1/2})$.

*Proof of part* (b) *for* $\hat{L}$. Now suppose that $f'$ is Hölder-continuous on $[a,b]$ with exponent $\beta - 1 \in (0,1]$ and constant $L > 0$. Here we define $U_n(x) := f(x) + \gamma\delta_n^\beta + 1\{x \notin [a+\delta_n, b-\delta_n]\}\infty$ with a constant $\gamma > 0$, and we assume that $\hat{U} \leq U_n$.

Suppose that $\hat{L}(t) \leq z := f(t) - 2\epsilon_n$ for some $t \in [a,b]$ and $\epsilon_n > 0$. If $t \leq b - 2\delta_n$, then

$$\begin{aligned}
h_{n,t,z}(t + \lambda\delta_n) &\leq z + \lambda(U_n(t + \delta_n) - z) \\
&= f(t) - 2\epsilon_n + \lambda(f(t + \delta_n) - f(t) + 2\epsilon_n + \gamma\delta_n^\beta) \\
&= f(t) - 2(1-\lambda)\epsilon_n + \lambda \int_0^{\delta_n} f'(t+s)\,ds + \lambda\gamma\delta_n^\beta
\end{aligned}$$

for $0 \leq \lambda \leq 1$. Thus

$$\begin{aligned}
&f(t + \lambda\delta_n) - h_{n,t,z}(t + \lambda\delta_n) \\
&= 2(1-\lambda)\epsilon_n + \lambda \int_0^{\delta_n} (f'(t + \lambda s) - f'(t + s))\,ds - \lambda\gamma\delta_n^\beta \\
&\geq 2(1-\lambda)\epsilon_n - \lambda \int_0^{\delta_n} L(1-\lambda)^{\beta-1} s^{\beta-1}\,ds - \lambda\gamma\delta_n^\beta \\
&\geq \epsilon_n - O(\delta_n^\beta)
\end{aligned}$$

uniformly for $0 \leq \lambda \leq 1/2$. Analogous arguments apply in the case $t \geq a + 2\delta_n$. Consequently there is an interval $[a_n, b_n] \subset [a,b]$ of length $\delta_n/2$ such that $f - h_{n,t,z} \geq \epsilon_n - O(\delta_n^\beta)$, provided that $a + 2\delta_n \leq b - 2\delta_n$. Again we choose $(d_n, j_n) \in \mathcal{T}_n$ with $d_n$ as large as possible such that $[x_{j_n - d_n + 1}, x_{j_n + d_n - 1}] \subset [a_n, b_n]$. Then $d_n \geq (c/4 + o(1))\delta_n n$, and (13) and (10) lead to (12). Thus $\epsilon_n = O(\delta_n^\beta) = O(\rho_n^{\beta/(2\beta+1)})$. □

## 6. Numerical examples

At first we illustrate the confidence band $(\hat{L}, \hat{U})$ defined in Section 3 with some simulated data. Precisely, we generated

$$Y_i = f(x_i) + \sigma\epsilon_i$$

with $x_i := (i - 1/2)/n$, $n = 500$ and

$$f(x) := \begin{cases} -12(x - 1/3) & \text{if } x \leq 1/3, \\ (27/2)(x - 1/3)^2 & \text{if } x \geq 1/3. \end{cases}$$

Moreover, $\sigma = 1/2$, and the random errors $\epsilon_1, \ldots, \epsilon_n$ have been simulated from a student distribution with five degrees of freedom. Figure 4 depicts these data together with the corresponding 95%–confidence band $(\hat{L}, \hat{U})$ and $f$ itself. Note that the width of the band is smallest near the center of the interval $(0, 1/3)$ on which $f$ is linear. This is in accordance with part (a) of Theorem 1.

Secondly we applied our procedure to a dataset containing the income $x_i$ and the expenditure $Y_i$ for food in the year 1973 for $n = 7125$ households in Great Britain (Family Expenditure Survey 1968–1983). This dataset has also been analyzed by Härdle and Marron [4]. They computed simultaneous confidence intervals for $\mathbb{E}(Y_i) = \tilde{f}(x_i)$ by means of kernel estimators and bootstrap methods. Figure 5 depicts the data. In order to enhance the main portion, the axes have been chosen such that 72 outlying observations are excluded from the display. Figure 6 shows a 95%–confidence band for the *isotonic* median function $f$, as described by Dümbgen and Johns [2]. Figure 7 shows a 95%–confidence band for the *concave* median function $f$, as described in the present paper. Note that the latter band has substantially smaller width than the former one. This is in accordance with our theoretical results about rates of convergence.

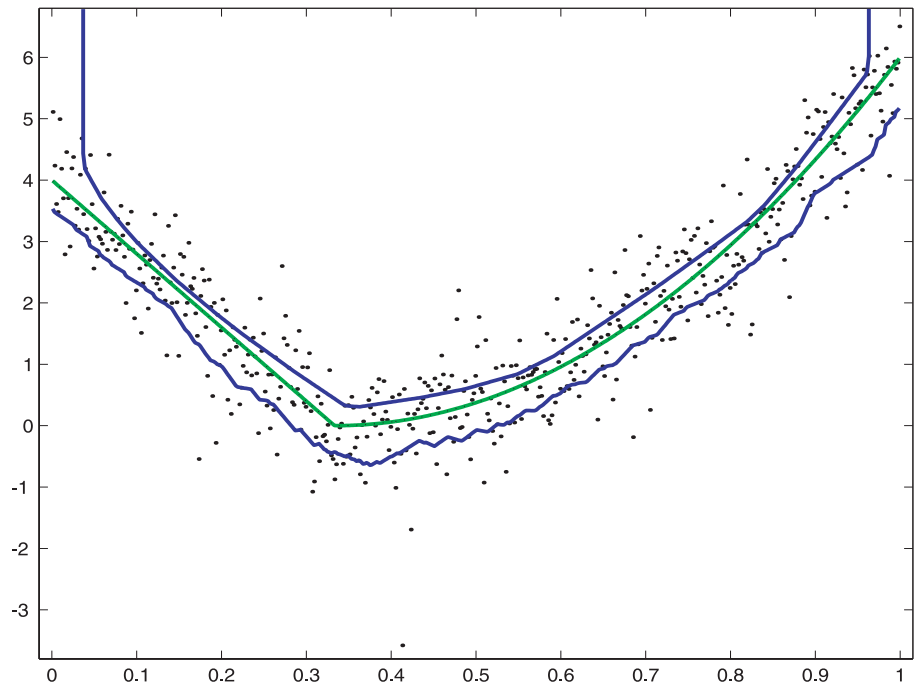

FIG 4. *Simulated data and 95%–confidence band $(\hat{L}, \hat{U})$, where $n = 500$.*

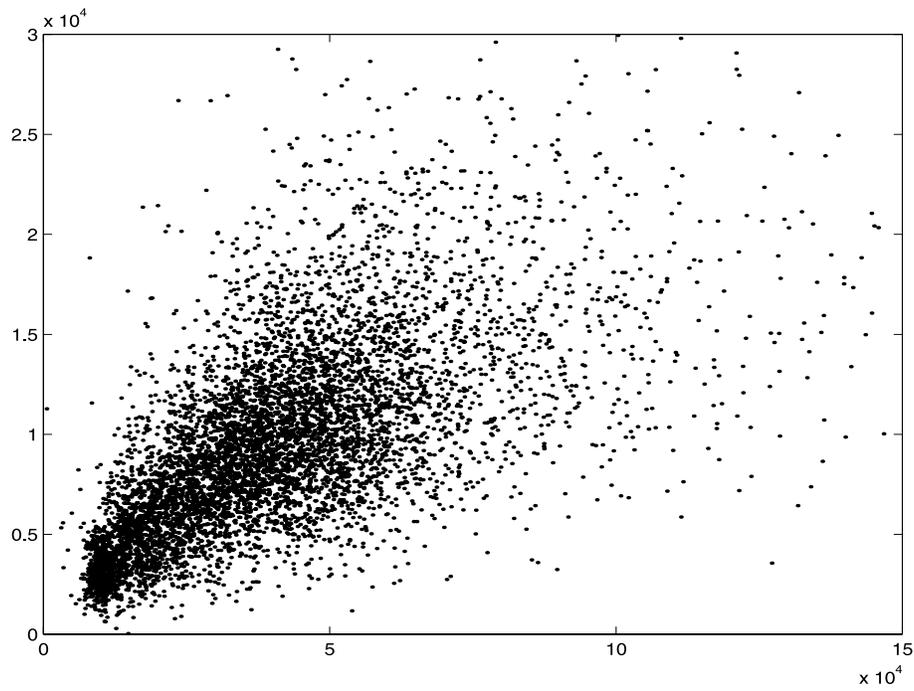

FIG 5. *Income-expenditure data.*

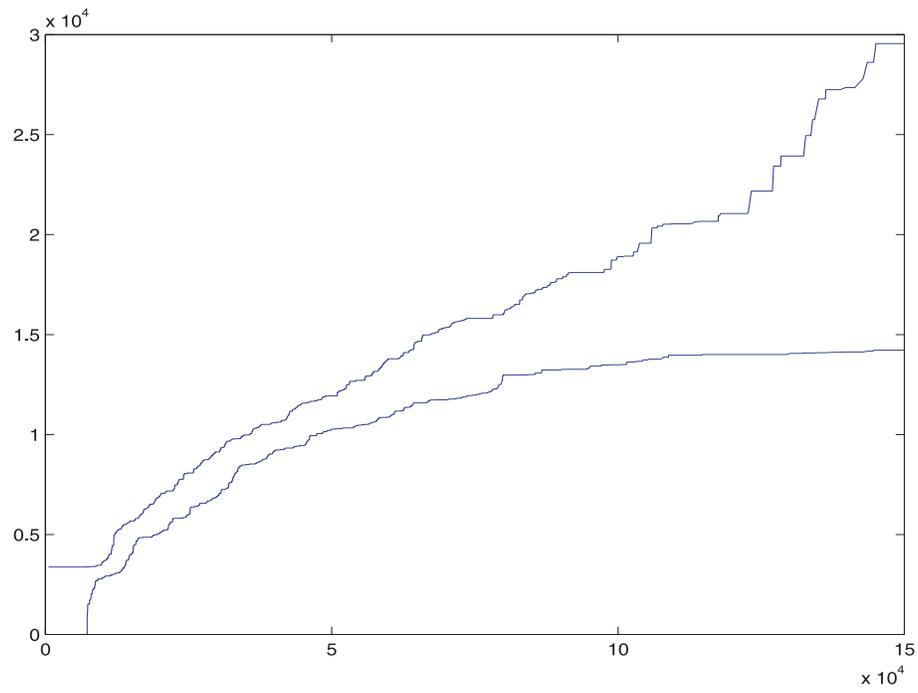

Fig 6. 95%–confidence band for isotonic median function.

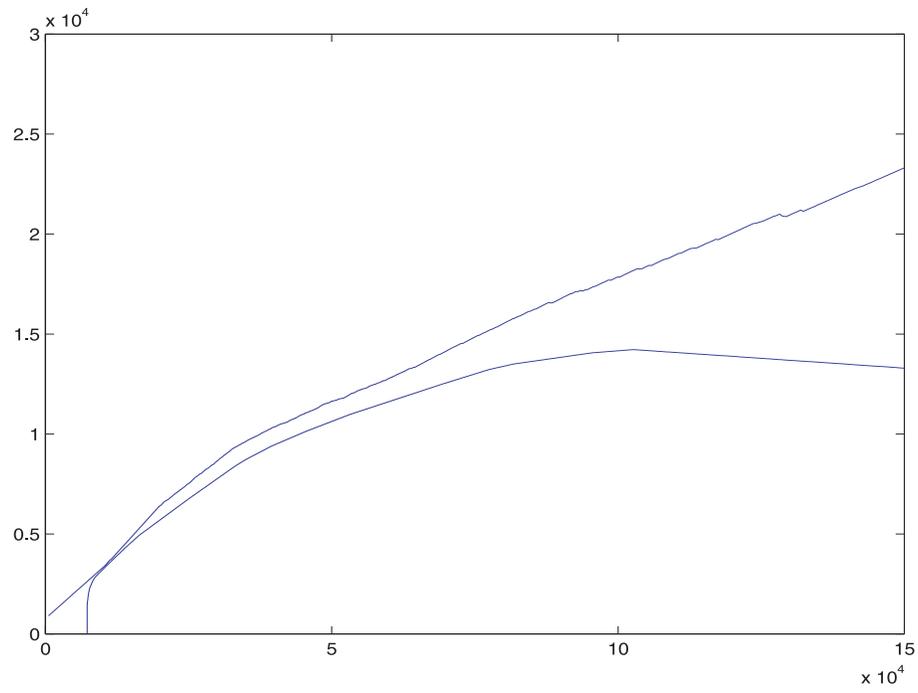

Fig 7. 95%–confidence band for concave median function.

**Acknowledgments.** The author is indebted to Geurt Jongbloed for constructive comments. Many thanks also to Wolfgang Härdle (Humboldt Unversity Berlin) for providing the family expenditure data.